\documentclass[conf]{new-aiaa} 

\usepackage{amsmath} 
\usepackage{amsfonts}
\newtheorem{definition}{Definition}[section]

\usepackage{float} 
\usepackage[skip=1pt,font=footnotesize]{caption}

\usepackage{subcaption}

\usepackage{etaremune}

\usepackage{tikz}
\usetikzlibrary{shapes,arrows,calc,positioning}
\tikzstyle{bigblock} = [draw, fill=blue!20, rectangle, 
    minimum height=6em, minimum width=8em]
\tikzstyle{medblock} = [draw, fill=blue!20, rectangle, 
    minimum height=4em, minimum width=4em]    
\tikzstyle{mux} = [draw, fill=black!20, rectangle, 
    minimum height=5em, minimum width=0.1em]    
\tikzstyle{smallblock} = [draw, fill=blue!20, rectangle, 
    minimum height=2em, minimum width=3em]
    
\tikzstyle{data_block} = [draw, fill=green!20, rectangle, 
    minimum height=2em, minimum width=3em]
\tikzstyle{ops_block} = [draw, fill=blue!20, rectangle, 
    minimum height=2em, minimum width=3em]    
\tikzstyle{est_block} = [draw, fill=red!20, rectangle, 
    minimum height=2em, minimum width=3em]    
    
\tikzstyle{sum} = [draw, fill=blue!20, circle, node distance=1cm,minimum height=0.5cm]
\tikzstyle{signal} = [coordinate]
\tikzstyle{pinstyle} = [pin edge={to-,thin,black}]
\tikzstyle{block} = [draw, fill=blue!20, rectangle, 
    minimum height=3em, minimum width=9em]
\tikzstyle{blockS} = [draw, fill=blue!20, rectangle, 
    minimum height=3em, minimum width=4em]    
\tikzstyle{input} = [coordinate]
\tikzstyle{output} = [coordinate]
\usetikzlibrary{matrix}

\usetikzlibrary{positioning}
\usetikzlibrary{math}

\usepackage{hyperref}
\usepackage{xcolor}
\hypersetup{
    colorlinks,
    linkcolor={blue!100!black},
    citecolor={blue!50!black},
    urlcolor={blue!80!black}
}

\newcommand{\bc}{\begin{center}}
\newcommand{\ec}{\end{center}}
\newcommand{\benum}{\begin{enumerate}}
\newcommand{\eenum}{\end{enumerate}}
\newcommand{\nn}{\nonumber}
\newcommand{\matl}{\left[ \begin{array}}
\newcommand{\matr}{\end{array} \right]}
\newcommand{\matls}{\left[ \begin{smallmatrix}}
\newcommand{\matrs}{\end{smallmatrix} \right]}
\newcommand{\isdef}{\stackrel{\triangle}{=}}

\newcommand{\pder}[2]{\mbox{$\dfrac{\partial {#1}}{\partial {#2}}$}}

\newcommand{\rmT}{{\rm T}}

\newcommand{\rmc}{{\rm c}}
\newcommand{\rmd}{{\rm d}}
\newcommand{\rme}{{\rm e}}

\newcommand{\rmm}{{\rm m}}

\newcommand{\BBR}{{\mathbb R}}

\newcommand{\SR}{{\mathcal R}}

\renewcommand{\matl}{\begin{bmatrix}}
\renewcommand{\matr}{\end{bmatrix} }

\usepackage{enumitem}
\newlist{todolist}{itemize}{2}
\setlist[todolist]{label=$\square$}
\usepackage{pifont}

\usepackage{tikz}
\usetikzlibrary{calc,patterns,decorations.pathmorphing,decorations.markings}
\usetikzlibrary{shapes,arrows}
\usepackage{verbatim}

\tikzstyle{block} = [draw, fill=blue!20, rectangle, 
    minimum height=3em, minimum width=6em]
\tikzstyle{sum} = [draw, fill=blue!20, circle, node distance=1cm]
\tikzstyle{input} = [coordinate]
\tikzstyle{output} = [coordinate]
\tikzstyle{pinstyle} = [pin edge={to-,thin,black}]

\usetikzlibrary{positioning}
\usetikzlibrary{math}

\usepackage{tikz}
\usetikzlibrary{shapes,arrows,calc,positioning}
\tikzstyle{bigblock} = [draw, fill=blue!20, rectangle, 
    minimum height=6em, minimum width=8em]
\tikzstyle{medblock} = [draw, fill=blue!20, rectangle, 
    minimum height=4em, minimum width=4em]    
\tikzstyle{mux} = [draw, fill=black!20, rectangle, 
    minimum height=5em, minimum width=0.1em]    
\tikzstyle{smallblock} = [draw, fill=blue!20, rectangle, 
    minimum height=3em, minimum width=4em]
\tikzstyle{sum} = [draw, fill=blue!20, circle, node distance=1cm]
\tikzstyle{signal} = [coordinate]
\tikzstyle{pinstyle} = [pin edge={to-,thin,black}]
\tikzstyle{block} = [draw, fill=blue!20, rectangle, 
    minimum height=3em, minimum width=6em]
\tikzstyle{blockS} = [draw, fill=blue!20, rectangle, 
    minimum height=3em, minimum width=4em]    
\tikzstyle{input} = [coordinate]
\tikzstyle{output} = [coordinate]
\usetikzlibrary{matrix}

\tikzset{add/.style n args={4}{
    minimum width=1mm,
    path picture={
        \draw[black, thick] 
            (path picture bounding box.south east) -- (path picture bounding box.north west)
            (path picture bounding box.south west) -- (path picture bounding box.north east);
        }
    }
}

\tikzset{Frame_into/.pic={
        code={\tikzset{scale=1}
        \tikzmath
            {
                \l  = 1;
                \Rc = 0.15;
            } 
        \draw [thick,->] (0, 0) -- +(\l, 0);
        \draw [thick,->] (0, 0) -- +(0, \l);
        \draw [thick, fill=white] 
			    (0,0) circle [radius=\Rc];
	    \draw [thick] ({\Rc*cos(45)}, {0\Rc*sin(45)}) 
	               -- ({\Rc*cos(225)}, {0\Rc*sin(225)});
        \draw [thick] ({\Rc*cos(135)}, {0\Rc*sin(135)}) 
	               -- ({\Rc*cos(315)}, {0\Rc*sin(315)});
  }}
}

\tikzset{Frame_outof/.pic={
        code={\tikzset{scale=1}
        \tikzmath
            {
                \l  = 1;
                \Rc = 0.15;
            } 
        \draw [thick,->] (0, 0) -- +(\l, 0);
        \draw [thick,->] (0, 0) -- +(0, \l);
        \draw [thick, fill=white] 
			    (0,0) circle [radius=\Rc];
	    \draw [thick, fill=black] 
			    (0,0) circle [radius={0.2*\Rc}];
  }}
}

\usepackage[style=ieee ,sorting=none,maxbibnames=99,giveninits, doi=false]{biblatex}
\title{\LARGE \bf
Longitudinal Flight Dynamics Control Based on Feedback Linearization and Normal Canonical Form
}

\title{\LARGE \bf
MIMO Input-Output Linearization with Applications for Longitudinal Flight Dynamics
}

\title{\LARGE \bf
Input-Output Linearization of MIMO Systems \\ with Application to Longitudinal Flight Dynamics
}

\title{\LARGE \bf
Longitudinal Flight Dynamics Control based on Input-Output Linearization despite Zero Dynamics
}

\title{\LARGE \bf
Circumventing Unstable Zero Dynamics in \\ Input-Output Linearization-based Control of \\ Longitudinal Flight Dynamics
}

\title{\LARGE \bf
Circumventing Unstable Zero Dynamics in \\ Input-Output Linearization of Longitudinal Flight Dynamics
}

\author{Jhon Manuel Portella Delgado and Ankit Goel
\thanks{Jhon Manuel Portella Delgado is a graduate student in the Department of Mechanical Engineering, University of Maryland, Baltimore County, 1000 Hilltop Circle, Baltimore, MD 21250. {\tt\small jportel1@umbc.edu}}%
\thanks{Ankit Goel is an Assistant Professor in the Department of Mechanical Engineering, University of Maryland, Baltimore County, 1000 Hilltop Circle, Baltimore, MD 21250. {\tt\small ankgoel@umbc.edu }}%
}
\usepackage{graphicx}      
\usepackage{textcomp}
\usepackage{calc}
\usepackage{mathtools}
\usepackage{xparse}
\usepackage{amsmath}
\usepackage{cases}
\usepackage{mathtools}
\usepackage{cuted}

\begin{document}

\maketitle
\begin{abstract}
In this paper, we consider the problem of input-output linearization of the longitudinal flight dynamics. 
In longitudinal flight dynamics, inputs are typically thrust and elevator deflection whereas the outputs are the velocity and the flight path angle. 
An input-output linearization-based controller can be designed to render the multi-input, multi-output system linear; however, the resulting zero dynamics turns out to be unstable.
In this work, we remove the zero dynamics from the closed-loop dynamics by considering an additional output. 
Although the additional output makes the system tall, which, in general, means that the input-to-output dynamics can not be linearized, we show that in the case of longitudinal flight dynamics, linearization is possible due to special geometric properties of the nonlinear terms. 


%


\end{abstract}

\section{INTRODUCTION}

Although input-output linearization (IOL) methods have been well-studied for square systems and nonsquare systems with fewer outputs than inputs, these methods have not been explored for systems with more outputs than inputs \cite{kolavennu2001nonlinear,alharbi2019backstepping,Khalil:1173048,isidori1985nonlinear}.  
In the classical input-output linearization method, a diffeomorphism is used to transform the state such that the dynamics matrix of a part of the transformed state is in the Jordan form, whereas the input matrix potentially remains a nonlinear function of the full state \cite{Khalil:1173048}. 
It turns out that, in the case of square systems and nonsquare systems with fewer outputs than inputs, the transformed input matrix, which is square or wide, is usually full-column rank. 
Nonlinearities thus can be canceled exactly to yield linear input-output dynamics.
However, the resulting zero dynamics may be unstable, rendering the IOL-based control impractical.  

The undesirable zero dynamics may be circumvented by considering additional outputs. 
However, additional outputs may make the system tall.  
In the case of nonsquare systems with more outputs than inputs, the transformed input matrix, in general, can not be inverted to cancel the nonlinearities, again  rendering the IOL-based control impractical.  
In this paper, we show that IOL method can, however, be applied to linearize a class of tall nonlinear systems. 


The work presented in this paper is motivated by the problem of linearizing the longitudinal aircraft dynamics. 
The longitudinal dynamics of an aircraft typically has two inputs, namely thrust and the elevator-deflection angle and two outputs, namely the velocity and the flight-path angle \cite{etkin1996dynamics,stevens2003aircraft}.
Although this square system can be linearized by input-output linearization methods, the resulting zero dynamics turns out to be unstable. 
By considering an additional output in the linearization process, we show that zero dynamics can be eliminated in this case.
%
%

Several methods have been explored to regulate the states of aircraft longitudinal dynamics.
The total energy control system proposed in \cite{lambregts1983integrated,lambregts1983vertical} transforms the states to energy states and uses heuristically tuned PID gains to regulate the states of the aircraft. 
However, this approach does not guarantee the stability of the closed-loop system and does not provide a straightforward mechanism to tune the transient response.
Nonlinear backstepping methods have also been investigated to solve this problem \cite{gavilan2011control} to obtain stability guarantees. 
However, since backstepping methods require the dynamics to be in a strict feedback form, these approaches often omit the effect of the elevator deflection on the lift in order to formulate the dynamics in the strict feedback form.
This paper considers a more realistic aircraft dynamics by including the effect of elevator deflection on the lift.
Since the dynamics considered in this paper is not in a strict feedback form, classical backstepping methods are not applicable.

The contribution of this paper is thus the extension of the input-output linearization method to a class of nonsquare systems with more outputs than inputs, its application to design a controller to linearize the aircraft longitudinal dynamics, and numerical demonstration of the proposed method in a nominal and off-nominal scenario.  
%
%
%
The paper is organized as follows. 
Section \ref{sec:LAD} describes the longitudinal aircraft dynamics used to design the input-output linearizing controller,  
Section \ref{sec:MIMO_IOL} presents the input-output linearization method for MIMO systems,
Section \ref{sec:IOL2LAD} shows the application of the input-output linearization method to the problem of linearizing longitudinal aircraft dynamics in the case of two and three outputs, and 
Section \ref{sec:simulations} shows the results of the numerical simulations of the closed-loop longitudinal aircraft dynamics.
Finally, the paper concludes with a discussion of results and future research directions in Section \ref{sec:conclusions}.





\section{Longitudinal Aircraft Dynamics}
\label{sec:LAD}
This section reviews the longitudinal dynamics of an aircraft and presents the notation used in this paper. 
The longitudinal flight dynamics are given by
\begin{align}
        \dot{V}&=\frac{1}{m}[F \cos(\alpha) - D - mg \sin(\gamma)]
        \label{eq:Vdot}
        ,
        \\
        \dot{\gamma}&=\frac{1}{mV}[F \sin(\alpha) + L -mg \cos(\gamma)]
        ,
        \\
        \dot{\theta}&=q
        \\
        \dot{q}&=\frac{M}{I_{yy}},
        \label{eq:qdot}
    \end{align}
where 
$V$ is the velocity, 
$\gamma$ is the flight-path angle,
$\theta$ is the pitch angle,
$\alpha \isdef \theta-\gamma$ is the angle-of-attack,
$q$ is the pitch rate, 
$F$ is the thrust, and 
$\delta_\rme$ is the elevator deflection angle 
\cite{anandakumar2022adaptive}. 
The lift $L,$ the drag $D,$ and the moment $M$ are parameterized as
\begin{align}\label{eq:LDM_def}
    L&=\frac{1}{2}\rho V^2 S C_\ell,\quad
    D =\frac{1}{2}\rho V^2 S C_\rmd,\quad
    M =\frac{1}{2}\rho V^2 S \overline{c} C_\rmm,
\end{align}
where $\rho$ is the air density, 
$S$ is the wing surface area, and 
$\overline{c}$ is the mean chord length.
Finally, the lift coefficient $C_\ell,$ the drag coefficient $C_\rmd,$ and the moment coefficient $C_\rmm$ are parameterized as
\begin{align}\label{eq5}
    C_\ell&= C_{\ell,0} + C_{\ell,\alpha}~\alpha + C_{\ell,{\delta_\rme}}~\delta_e,\\
    C_\rmd&= C_{\rmd,0} + C_{\rmd,\alpha}~\alpha,\\
    C_\rmm&= C_{\rmm,0} + C_{\rmm,\alpha}~\alpha + C_{\rmm,\delta_{\rm{e}}}~\delta_ {\rm{e}},
\end{align}
where 
$ C_{\ell,0}, $
$C_{\ell,\alpha}, $
$C_{\ell,{\delta_{\rm{e}}}},$ 
$C_{\rmd,0}, $
$C_{\rmd,\alpha}, $
$ C_{\rmm,0}, $ 
$ C_{\rmm,\alpha}$, and $C_{\rmm,\delta_{\rm{e}}}$ are aircraft aerodynamic coefficients, and $\delta_{\rm{e}}$ is the elevator angle.
%
Note that the inclusion of $C_{\ell,{\delta_\rme}}$ in the dynamics makes the system non-triangular, and hence, conventional backstepping methods can not be applied to this problem \cite{gavilan2011control}.

\begin{figure}[!ht]
    \centering
    \includegraphics[width = 0.6\columnwidth]{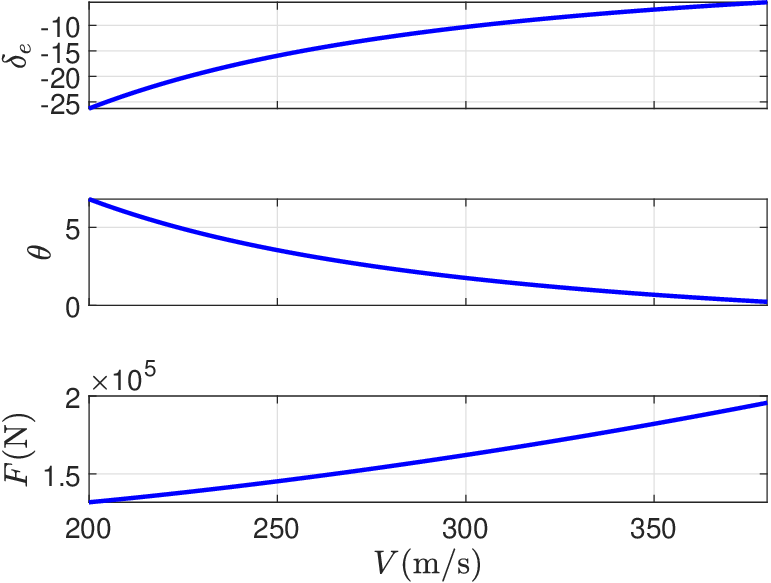}
    \caption{Trim conditions at various aircraft velocities. Note that the elevator deflection and the pitch angles are expressed in degrees. }
    \label{fig:trim_condition}
\end{figure}

%
%
In this work, we consider the physical parameters of the A330 aircraft given in \cite{abdulhamitbilal2007matlab} to simulate the longitudinal flight dynamics.
The various parameters are shown in Table \ref{table_parameters}.
\begin{table}[!ht]
\centering
\begin{tabular}{|l|l|l|} 
\hline
Parameter & value & Description                   \\ 
\hline
$C_{\ell,0}$      & 0.2301 & Lift coefficient                  \\
$C_{\ell,\alpha}$       & 5.9598 & Lift coefficient                 \\
$C_{\ell,\delta_e}$      & 0.2391 & Lift coefficient \\
$C_{\rm{m,0}}$       & -0.0812  & Pitching moment coefficient               \\
$C_{\rm{m,\alpha}}$      & -3.1069  & Pitching moment coefficient\\
$C_{\rm{m,\delta_e}}$     & -0.9816   & Pitching moment coefficient              \\
$C_{\rm{d,0}}$      & 0.0172  & Drag coefficient\\
$C_{\rm{d,\alpha}}$      & 0.2223  & Drag coefficient               \\
$S$        & 363.12 $\rmm^2$  & Area of the wing               \\
$m$       & 254,842~$\rm{Kg}$    & Mass             \\
$\rho$      & $0.4127$ $\rm kg/m^3$  & Air density \\
$\overline{c}$       & 7.49 $\rmm$    & Chord length of the wing                \\
g        & $9.81~\rm{m/s^2}$    & Acceleration due to gravity                \\
$I_{yy}$      & $30,513,547 ~\rm{kgm^2}$ & Moment of inertia            \\
\hline
\end{tabular}
\caption{Physical parameters used to simulate  Longitudinal dynamics.
}
\label{table_parameters}
\end{table}
The trim conditions are obtained numerically by setting the derivative of the state to zero in the system's dynamics. 
Figure \ref{fig:trim_condition} shows the trim conditions at a steady state for several velocities.
%

\section{MIMO Input-Output Linearization}
\label{sec:MIMO_IOL}
This section reviews the multi-input, multi-output extension of the input-output linearizing control presented in \cite{kolavennu2001nonlinear}, and extends it to the case of nonsquare systems with more outputs than inputs.
Consider an affine system
\begin{align}
    \dot x 
        &= 
            f(x) + g(x) u,
    \label{eq:xdot_gen}
    \\
    y 
        &=
            h(x),
    \label{eq:y_gen}
\end{align}
where $x(t)\in \BBR^{l_x}$ is the state, 
$u(t)\in \BBR^{l_u}$ is the input, 
$y(t)\in \BBR^{l_y}$ is the output,
and $f, g, h$ are smooth functions of appropriate dimensions. 
The objective is to construct a control law $u=k(x)+v,$ where $v$ is the intermediate control, such the dynamics from the intermediate input $v$ to the output is linear, that is, 
\begin{align}
    \dot \xi &= A \xi + B v, \\ 
    y   &= C \xi, 
\end{align}
where $A,B,C$ are desired matrices. 
Finally, the intermediate control $v$ can be designed to obtain the desired output response using tools from linear systems theory. 

The following definitions appear in \cite{isidori1985nonlinear} and are repeated here for further use in the paper. 


\begin{definition}
    \label{def:rel_degree}
    In the system \eqref{eq:xdot_gen}, \eqref{eq:y_gen}.
    the relative degree of the $i$th output $y_i$ is the smallest integer $\rho_i \ge 0$ such that $\rho_i$-th derivative of $y_i,$ that is
    $y_i^{(\rho_i)},$ is an explicit function of input $u$. 
\end{definition}

\begin{definition}
    The relative degree of the system \eqref{eq:xdot_gen}, \eqref{eq:y_gen} is the sum of the relative degree of each of its outputs, that is, $\rho \isdef \sum_i^{l_y} \rho_i.$
\end{definition}

\begin{definition}
    Let $\zeta \colon \BBR^{l_x} \to \BBR^{l_\zeta}$
    and 
    $f \colon \BBR^{l_x} \to \BBR^{l_f}$
    be smooth functions. Then, the Lie derivative of $\zeta$ with respect to $f$, denoted by $L_f \zeta,$ is
    \begin{align}
        L_f \zeta(x)
            \isdef
                \pder{\zeta(x)}{x} ^\rmT f(x).
    \end{align}
    Note that $L_f \zeta \colon \BBR^{l_x} \to \BBR^{l_\zeta}.$
\end{definition}



%

Consider the transformation
\begin{align}\label{eq_T}
    T : \BBR^{l_x} &\to \BBR^{l_x}
    \nn \\
        T(x)
            &=
                \matl
                    \phi(x) \\
                    \psi(x)
                \matr,
\end{align}
where $\phi(x) $ satisfies 
\begin{align}\label{eq_phi}
    L_g \phi(x) = 0,
\end{align}
and
\begin{align}
    \psi(x) 
        =
            \matl
                \psi_1(x) \\
                \vdots \\
                \psi_{l_y}(x)
            \matr,
    \label{eq:psi_def}
\end{align}
where
\begin{align}
    \psi_i(x)
        \isdef
            \matl
                h_i(x) \\
                L_f h_i(x) \\
                \vdots \\
                L_f^{\rho_i-1} h_i(x)
            \matr
        \in 
            \BBR^{\rho_i}.
    \label{eq:psi_i_def}
\end{align}
Note that $\phi : \BBR^{l_x} \to \BBR^{{l_x} - \rho}$ and, for $i=1, \ldots, l_y,$
$\psi_i : \BBR^{l_x} \to \BBR^ {\rho_i},$ and thus
$\psi : \BBR^{l_x} \to \BBR^ \rho.$
Furthermore, the functions $\psi_i$ are well-defined since the functions $f,g,h$ are assumed to be smooth. 
However, $\phi$ satisfying \eqref{eq_phi} may or may not exist. 
%
Assuming that $\phi$ satisfying \eqref{eq_phi} exists and defining $\eta \isdef \phi(x), $ it follows 
that 
\begin{align}
    \dot \eta 
        &=
            L_f \phi(x) + L_g \phi(x) u
        =
            L_f \phi(x),
    \label{eq_eta_dot}
\end{align}
where $L_g \phi(x)=0$ by construction. 
Note that \eqref{eq_eta_dot} is the \textit{zero dynamics} \cite{Khalil:1173048} of the system \eqref{eq:xdot_gen}, \eqref{eq:y_gen}.

Next, defining $\xi \isdef \psi(x),$ it follows
that 
\begin{align}
    \dot \xi
        &=
            L_f \psi(x) + L_g \psi(x) u.
    \label{eq:xi_dot}
\end{align}
Next, note that 
\begin{align}
    L_f \psi(x)
        =
            A_\rmc \xi 
            +
            B_\rmc 
            \matl 
                L_f^{\rho_1} h_1(x) \\
                \vdots \\
                L_f^{\rho_{l_y}} h_{l_y}(x)
            \matr,
    \label{eq:Lf_psi}
\end{align}
where $A_\rmc = {\rm diag  } (A_{\rmc,1}, \ldots, A_{\rmc, l_y}) \in \BBR^{\rho \times \rho }$ and 
$B_\rmc = {\rm diag} (b_{\rmc,1}, \ldots, b_{\rmc,l_y}) \in \BBR^{\rho \times l_y} $ and, for $i=1,\ldots, l_y,$ 
\begin{align}
    A_{\rmc,i}
        &\isdef 
            \matl 
                0 & 1 & 0 & \cdots & 0  \\
                0 & 0 & 1 & \cdots & 0  \\
                \vdots & \vdots & \ddots & \ddots & \vdots \\
                0 & \ldots & \ldots & 0 & 1  \\
                0 & \ldots & \ldots & 0 & 0  \\
            \matr
        \in \BBR^{\rho_i \times \rho_i},
    \quad
    b_{\rmc,i}
        \isdef 
            \matl 
                0\\
                \vdots \\
                1
            \matr
        \in \BBR^{\rho_i}.
\end{align}
Furthermore, 
\begin{align}
    L_g \psi(x)
        &=
            B_\rmc  
            \matl 
                L_g L_f^{\rho_1-1} h_1(x) \\
                \vdots \\
                L_g L_f^{\rho_{l_y}-1} h_{l_y} (x)
            \matr.
    \label{eq:Lg_psi}
\end{align}
Substituting \eqref{eq:Lf_psi} and \eqref{eq:Lg_psi} in \eqref{eq:xi_dot} thus yields
\begin{align}
    \dot \xi
        &=
            A_\rmc \xi 
            +
            B_\rmc
            \left(
                \alpha(x)
                +
                \beta(x) u
            \right),
\end{align}
where 
\begin{align}
    \alpha(x)
        &\isdef
            \matl 
                L_f^{\rho_1} h_1(x) \\
                \vdots \\
                L_f^{\rho_{l_y}} h_{l_y}(x)
            \matr
            \in \BBR^{l_y}, 
    \quad 
    \beta(x)
        \isdef 
            \matl 
                L_g L_f^{\rho_1-1} h_1(x) \\
                \vdots \\
                L_g L_f^{\rho_{l_y}-1} h_{l_y} (x)
            \matr
            \in \BBR^{l_y \times l_u}.
\end{align}
Finally, letting 
\begin{align}
    u(x) = -\beta(x)^{+}\alpha(x) + v,
    \label{eq:IOL_control_law}
\end{align}
yields
\begin{align}
    \dot \xi
        &=
            A_\rmc \xi 
            +
            B_\rmc  \Lambda(x) \alpha (x)
            +
            B_\rmc  v ,
    \label{eq:IOlinearizedsystem}
\end{align}
where
$\Lambda(x) \isdef I_{l_y} - \beta(x) \beta(x)^+ \in \BBR^{l_y \times l_y}$ and $v \in \BBR^{l_y}.$
Note that 
$\Lambda(x) 
        =
            {\rm diag  } (\lambda_1(x), \ldots, \lambda_{l_y}(x)),$
where, for $i = 1,\ldots, l_y,$ $\lambda_i(x) $ is either 1 or 0.

\subsection{Square and Wide Plants}
In square and wide plants, that is, $l_u \ge l_y$,
if, for all $x \in \BBR^{l_x}$, $\beta(x)$ is full-column rank, then, for all $x \in \BBR^{l_x}$, $\Lambda (x) = 0.$ 
In this case, 
\eqref{eq:IOL_control_law} is the \textit{input-output linearizing} (IOL) controller
and \eqref{eq:IOlinearizedsystem} is the \textit{input-output linearized system}.  
%
Consequently, all outputs can be directly manipulated by appropriately defining the intermediate control $v.$

\subsection{Tall Plants}
In tall plants, that is, $l_u < l_y$, 
$\Lambda(x) \neq 0.$
In fact, at least $l_y-l_u$ diagonal elements of $\Lambda(x)$ are equal to one at each instant. 
Furthermore, if $\beta(x)$ is full-column rank for all $x$, then exactly $l_y-l_u$ diagonal elements of $\Lambda(x)$ are equal to one. 
%
%
Finally, in the case where, for all $x\in \BBR^{l_x},$  $\alpha(x) \in \SR(\beta(x))$,  $\Lambda(x) \alpha(x) = 0.$ 
In this special case, the tall system is thus \textit{input-output linearizable}.

\section{Input-Output Linearization of Longitudinal Aircraft Dynamics}
\label{sec:IOL2LAD}
This section applies the input-output linearizing control presented in Section \ref{sec:MIMO_IOL} to the  problem of linearization of longitudinal aircraft dynamics. 
%
%
%
We first show that in the case of two outputs, the zero dynamics associated with the linearization is unstable. 
%
%
%
To circumvent the unstable zero dynamics, we consider an additional output, namely, the pitch angle.
With the additional output, we show that there is no zero dynamics.
Note that the pitch angle reference is required to implement this controller, which can be obtained using the trim computations, but in general, is not precisely known.  
However, as shown in the numerical example considered in the paper, an incorrect pitch angle reference in the controller yields nonzero but bounded steady-state errors.

%


\subsection{Two outputs}

Defining $x \isdef \matl V - \overline V & \gamma - \overline \gamma & \theta & q \matr^\rmT, $
where $\overline V$ and $\overline \gamma$ are the reference velocity and the flight path angle to be tracked, 
and 
    $u \isdef 
        \matl
            F & 
            \delta_{\rm{e}}
        \matr^\rmT,$
it follows that \eqref{eq:Vdot}-\eqref{eq:qdot} can be written as \eqref{eq:xdot_gen}, \eqref{eq:y_gen},
where 
\begin{align}
    f_1(x) &\isdef
            -\dfrac{\rho (x_1+\overline{V})^2 S}{2m}(C_{\rmd,0}+C_{\rmd,\alpha}(x_3-x_2-\overline{\gamma}))
            -            g\sin{(x_2+\overline {\gamma})}-\dot{\overline{V}} ,
            \label{eq_f1_2}
                \\
    f_2(x) &\isdef 
            \dfrac{\rho (x_1+\overline{V}) S}{2m}(C_{\ell,0}+C_{\ell,\alpha}(x_3-x_2-\overline{\gamma})) 
                - \dfrac{g\cos{(x_2+\overline{\gamma})}}{(x_1+\overline{V})} - \dot{\overline{\gamma}},
                \label{eq_f2_2}
                \\
    f_3(x) &\isdef x_4,
                \label{eq_f3_2}
                \\
    f_4(x) &\isdef 
                \dfrac{\rho (x_1+\overline{V}) S \overline{c}}{2I_{yy}}(C_{\rmm,0}+C_{\rmm,\alpha}(x_3-x_2-\overline{\gamma})),
                \label{eq_f4_2}
\end{align}
and
\begin{align}\label{eq_gx_sys}
g(x)=
    \begin{bmatrix}
        \dfrac{\cos{(x_3-x_2-\overline
        \gamma)}}{m} 
        & 0
            \\
        \dfrac{\sin{(x_3-x_2-\overline \gamma)}}{m(x_1+\overline {V})}
        &
        \dfrac{\rho (x_1+\overline{ V})S}{2m}C_{\ell,{\delta_{\rm{e}}}}
            \\
        0 & 0 \\
        0 & \dfrac{\rho (x_1+\overline {V}) S \overline{c}}{2I_{yy}}C_{\rmm,\delta_{\rm{e}}}
    \end{bmatrix}.
\end{align}
It is assumed that $\overline V, \dot {\overline{V}}, \overline \gamma, \dot {\overline \gamma}$ are well-defined. 
%
Consider the output 
\begin{align}
    y 
        =
            h(x)
        \isdef
            \matl 
                x_1 \\
                x_2 
            \matr.
    \label{eq:output1}
\end{align}
where $x_1$ is the velocity error and $x_2$ is the flight-path angle error.
Since the reference values of the velocity and the flight-path angle are generated by the guidance law, they are known. Along with the velocity and the flight-path angle measurements, $x_1$ and $x_2$ are thus known. 
Next, note that $\rho_1=\rho_2= 1,$ and thus $\rho=2,$
which implies that $\eta = \phi(x) \in \BBR^2$ and $\xi = \psi(x) \in \BBR^2.$
It follows from \eqref{eq:psi_i_def} that 
\begin{align}
    \xi_1
    &=
        x_1,
        \label{psi_1_2states}
        \\
    \xi_2
    &=
        x_2.
\end{align}
Furthermore, 
\begin{align}
    \beta(x)
        =
        \matl
        \dfrac{\cos{(x_3-x_2-\overline \gamma)}}{m} & 0\\
        \dfrac{\sin{(x_3-x_2-\overline \gamma)}}{m(x_1+\overline V)} & \dfrac{\rho (x_1+\overline V) S }{2m}C_{\ell,\delta_{\rm{e}}}
        \matr.
\end{align}
Note that, if $x_1 \neq -\overline{V}$ and $x_3-x_2-\overline{\gamma} \neq \dfrac{\pi}{2},$ then ${\rm det}(\beta(x)) \neq 0$.
The IOL control law \eqref{eq:IOL_control_law} thus yields 
\begin{align}
    \dot x_1 &= v_1, \\
    \dot x_2 &= v_2,
\end{align}
which implies that $x_1$ and $x_2$ can be arbitrarily regulated.

Next, solving \eqref{eq_phi} yields
\begin{align}
    \eta_1 &=
        \sin{(x_3-x_2-\overline{\gamma})}(x_1+\overline{V}),
       \label{eq_phi_sol1} 
        \\
    \eta_2
    &=
        \dfrac{m\overline{c} C_{\rmm,\delta_{\rm{e}}}}{I_{yy}C_{\ell,{\delta_{\rm{e}}}}}x_2-x_4,
        \label{eq_phi_sol2} 
\end{align}
and thus the zero dynamics is given by 
%
\begin{align}
    \dot \eta_1 
        &
        =
            \eta_1\bigg[-\dfrac{\rho\overline{V}S}{2m}\bigg(C_{\rmd,0}+C_{\rmd,\alpha}\sin^{-1}{\bigg(\dfrac{\eta_1}{\overline{V}}\bigg)}\bigg)
            -
            \dfrac{g\sin{\overline{\gamma}}}{\overline{V}}-\dfrac{\dot{\overline{V}}}{\overline{V}}\bigg]
            +\sqrt{(\overline{V}^2-\eta_1^2)}\bigg[-\dfrac{\rho\overline{V}S}{2m}\bigg(C_{\ell,0}+C_{\ell,\alpha}\sin^{-1}{\bigg(\dfrac{\eta_1}{\overline{V}}\bigg)}\bigg)
            \nn \\ &\quad
            +\dfrac{g\cos{\overline{\gamma}}}{\overline{V}}+\dot{\overline{\gamma}}-\eta_2\bigg],
    \\
    \dot \eta_2 
        &
        =
           \dfrac{\rho\overline{V}S\overline{c}C_{\rmm,\delta_{\rm{e}}}}{2I_{yy}C_{\ell,{\delta_{\rm{e}}}}}\bigg(C_{\ell,0}+C_{\ell,\alpha}\sin^{-1}{\bigg(\dfrac{\eta_1}{\overline{V}}\bigg)}\bigg)
           -\dfrac{mg\overline{c}C_{\rmm,\delta_{\rm{e}}}\cos{\overline{\gamma}}}{I_{yy}\overline{V}C_{\ell,{\delta_{\rm{e}}}}}-\dfrac{m\overline{c}C_{\rmm,\delta_{\rm{e}}}\dot{\overline{\gamma}}}{I_{yy}C_{\ell,{\delta_{\rm{e}}}}}
           -\dfrac{\rho\overline{V}S\overline{c}}{2I_{yy}}\bigg(C_{\rmm,0}+C_{\rmm,\alpha}\sin^{-1}{\bigg(\dfrac{\eta_1}{\overline{V}}\bigg)}\bigg).
\end{align}
Several simulations with nominal values of the parameters confirm the well-known fact that the zero-dynamics, in this case, is unstable. 
Furthermore, since
\begin{align}
    x_3&=
        \sin^{-1}{\left(\dfrac{\eta_1}{\xi_1+\overline{V}}\right)} + \xi_2 + \overline{\gamma},
        \label{eq_x3_eta}
        \\
    x_4&=
        \dfrac{m\overline{c} C_{\rmm,\delta_{\rm{e}}}}{I_{yy}C_{\ell,{\delta_{\rm{e}}}}}\xi_2 - \eta_2,
        \label{eq_x4_eta}
\end{align}
the pitch angle and the pitch rate also diverge due to  unstable zero dynamics. 

\subsection{Three outputs}
In order to remove the unstable zero dynamics, we consider the pitch angle as an additional output. 
In particular, we assume that the pitch angle of the aircraft is commanded to a desired value, which is assumed to be given by the desired trim condition. 
Redefining the state $x \isdef \matl V - \overline V & \gamma - \overline \gamma & \theta-\overline \theta & q \matr^\rmT,$
it follows that \eqref{eq:Vdot}-\eqref{eq:qdot} can be written as \eqref{eq:xdot_gen}, \eqref{eq:y_gen},
where 
\begin{align}
    f_1(x) &\isdef
            -\dfrac{\rho (x_1+\overline{V})^2 S}{2m}(C_{\rmd,0}+C_{\rmd,\alpha}(x_3+\overline{\theta}-x_2-\overline{\gamma}))
            -            g\sin{(x_2+\overline {\gamma})}-\dot{\overline{V}} ,
            \label{eq_f1}
                \\
    f_2(x) &\isdef 
            \dfrac{\rho (x_1+\overline{V}) S}{2m}(C_{\ell,0}+C_{\ell,\alpha}(x_3+\overline{\theta}-x_2-\overline{\gamma})) 
                - \dfrac{g\cos{(x_2+\overline{\gamma})}}{(x_1+\overline{V})} - \dot{\overline{\gamma}},
                \label{eq_f2}
                \\
    f_3(x) &\isdef x_4,
                \\
    f_4(x) &\isdef 
                \dfrac{\rho (x_1+\overline{V}) S \overline{c}}{2I_{yy}}(C_{\rmm,0}+C_{\rmm,\alpha}(x_3+\overline{\theta}-x_2-\overline{\gamma}))
               -\ddot{\overline{\theta}}
                \label{eq_f4}.
\end{align}
Consider the output 
\begin{align}
    y=h(x)\isdef
    \matl
        x_1\\
        x_2\\
        x_3
    \matr,
\end{align}
where $x_1, x_2, $ and $x_3$ are the velocity error, the flight-path angle error, and the pitch-angle error, respectively.
Note that $\rho_1=\rho_2= 1,$ and $\rho_3=2,$ and thus $\rho=4$ and $\xi = \psi(x) \in \BBR^4$.
Furthermore, it follows from \eqref{eq:psi_def} and \eqref{eq:psi_i_def} that $\psi(x)=  x$, and thus there is no zero dynamics. 
Furthermore, 
\begin{align}
    \beta(x)
        =
        \matl
        \dfrac{\cos{(x_3-x_2-\overline \gamma)}}{m} & 0\\
        \dfrac{\sin{(x_3-x_2-\overline \gamma)}}{m(x_1+\overline V)} & \dfrac{\rho (x_1+\overline V) S }{2m}C_{\ell, \delta_{\rm{e}}}\\
        0 & \dfrac{\rho (x_1+\overline {V}) S \overline{c}}{2I_{yy}}C_{\rmm,\delta_{\rm{e}}}
        \matr.
\end{align}

%
Finally, the IOL control law \eqref{eq:IOL_control_law} thus yields
\begin{align}
    \dot x
        =
            A_\rmc x 
            +
            B_\rmc  \Lambda(x) \alpha (x)
            +
            B_\rmc  v             
    \label{eq:LFD_IOL}
\end{align}
where 
\begin{align}
    A_\rmc 
        =
            \matl 
                0 & 0 & 0 & 0\\
                0 & 0 & 0 & 0\\
                0 & 0 & 0 & 1\\
                0 & 0 & 0 & 0
            \matr,
    B_\rmc 
        =
            \matl 
                1 & 0 & 0\\
                0 & 1 & 0\\
                0 & 0 & 0\\
                0 & 0 & 1
            \matr.
\end{align}

In \eqref{eq:Vdot}-\eqref{eq:qdot}, it turns out that 
$\Lambda(x) \alpha(x)$ is close to zero for all times, and thus the closed-loop dynamics is approximately 
\begin{align}
    \dot x_1 &= v_1, \\
    \dot x_2 &= v_2, \\
    \dot x_3 &= v_3. 
\end{align}
Letting 
\begin{align}
    v = 
        \matl
              -K_1 x_1\\
              -K_2 x_2\\
              -K_3 x_3 -K_4x_4,
        \matr
\end{align}
with $K_1,K_2,K_3,K_4>0$ ensures that $x \to 0,$ that is, velocity error, flight path angle error, and the pitch angle error converge to zero asymptotically.

\section{Simulation Results}
\label{sec:simulations}
In this section, we apply the IOL controller, presented in Section \ref{sec:IOL2LAD},  to regulate an aircraft's velocity, flight path angle, and pitch angle. 
%
%
%
%
%
In practice, the exact pitch reference may not be known.
However, as shown in the numerical simulations, the output error remains bounded in the case of an unknown bias in the pitch reference. 

\begin{figure}[!ht]
    \centering
    \includegraphics[width = 0.6\columnwidth]{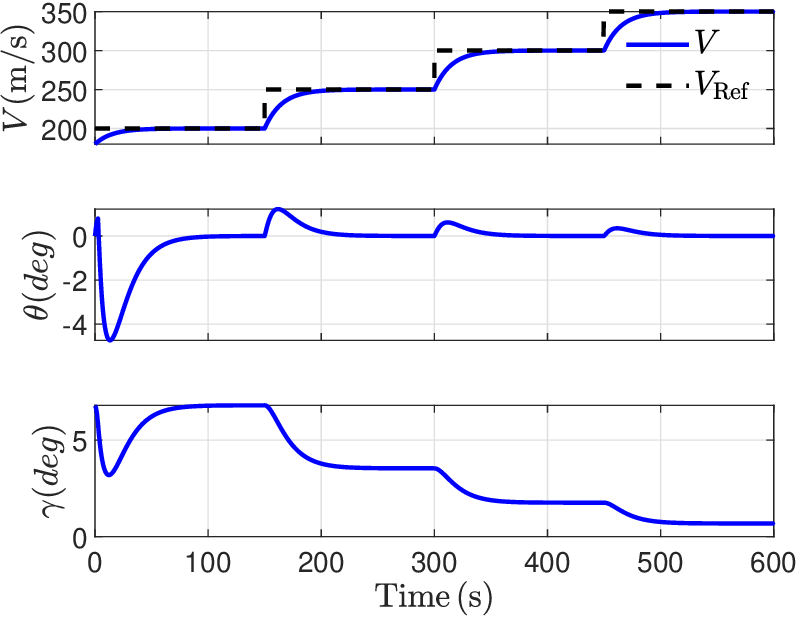}
    \caption{Velocity, flight-path angle, and the pitch angle response of the longitudinal aircraft dynamics with the linearizing controller. 
    Note that the output is shown in solid blue and the corresponding reference is shown in dashed black.}
    \label{fig:Signals}
\end{figure}
\begin{figure}[!ht]
    \centering
    \includegraphics[width = 0.6\columnwidth]{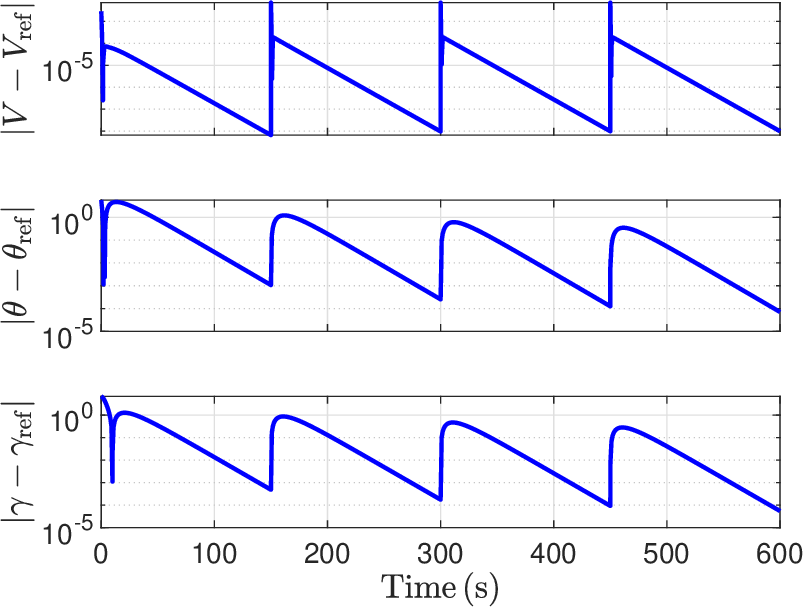}
    \caption{Absolute value of the velocity error, flight-path error, and the pitch angle error in the closed-loop simulation on a logarithmic scale. }
    \label{fig:Errors}
\end{figure}
\begin{figure}[!ht]
    \centering
    \includegraphics[width = 0.6\columnwidth]{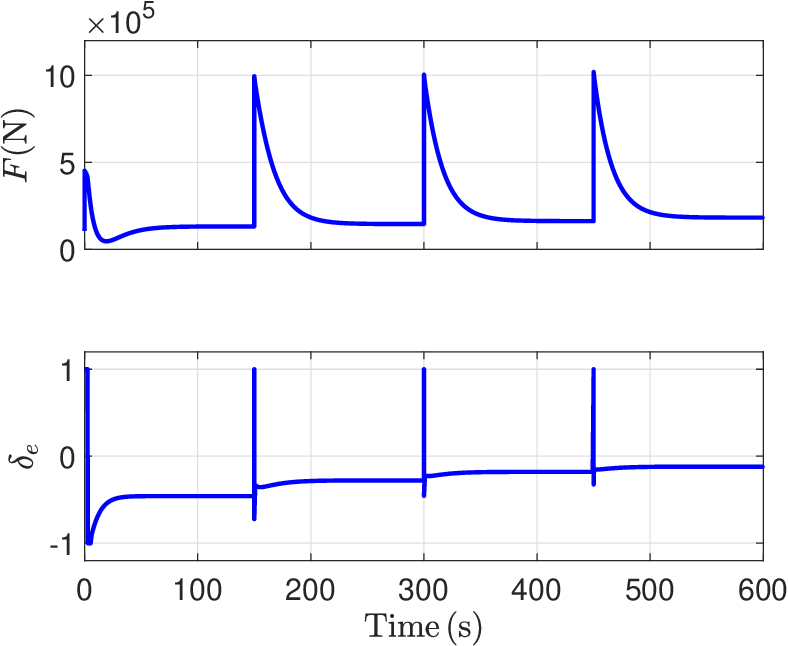}
    \caption{Thrust and elevator-deflection angle given by the linearizing controller.}
    \label{fig:Control_inputs_and_determinant}
\end{figure}

The aircraft is assumed to be flying at a steady state with a velocity of $180$ m/s at an altitude of $10,000$ km. 
The aircraft is then commanded to increase its velocity in steps every $150$ seconds. 
%
In the IOL controller, we set 
$k_1 = 4, k_2 = 1, k_3 = 30, k_4=200.$
Note that the IOL controller completely linearizes and decouples the dynamics of each output, thus the gains can be chosen to satisfy the desired transient requirements. 
%
Figure \ref{fig:Signals} shows the velocity, flight path angle, and pitch of the aircraft
and 
Figure \ref{fig:Errors} shows the absolute values of velocity error, flight path angle error, and pitch error on a logarithmic scale. 
Note that the errors go to zero exponentially. 
Figure \ref{fig:Control_inputs_and_determinant} shows the corresponding control inputs given by the IOL controller.

\begin{figure}[!ht]
    \centering
    \includegraphics[width = 0.6\columnwidth]{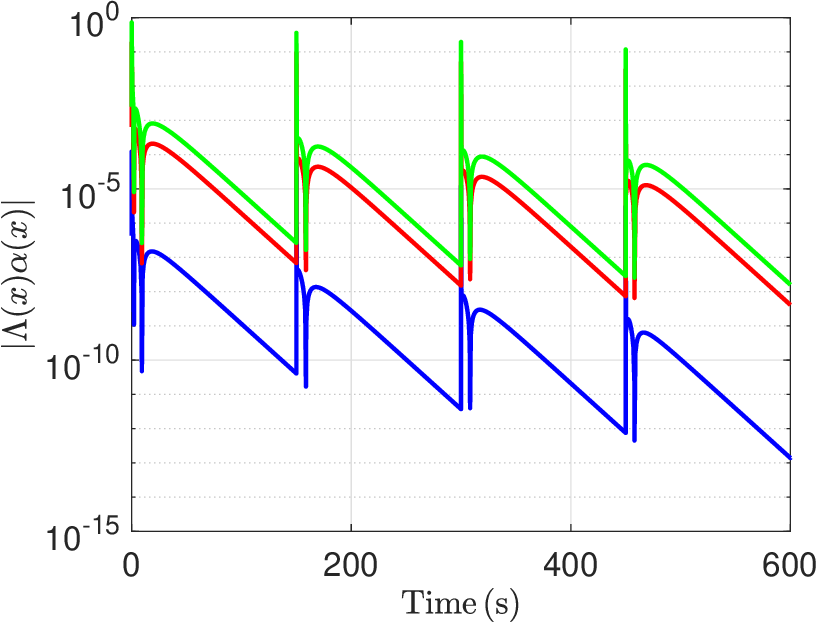}
    \caption{Components of $|\Lambda(x)\alpha(x)|$ on a logarithmic scale. Note that asymptotic convergence $\Lambda(x)\alpha(x)$ to $0$ allows complete linearization of the tall MIMO system. }
    \label{fig:LambdaGamma}
\end{figure}

Since $\beta(x) $ is a $3\times 2$ matrix, rank of $\Lambda(x)$ is at least one at all times, therefore, the nonlinear term in \eqref{eq:IOlinearizedsystem} is not necessarily zero. 
However, in this application, $\beta(x)$ is \textit{almost} in the range space of the columns of $\alpha(x),$ and thus $\Lambda(x) \alpha(x)$ is approximately zero at all times.  
Figure \ref{fig:LambdaGamma} shows the components of $\Lambda(x) \alpha(x)$ on a logarithmic scale. 
Note that since $\Lambda(x) \alpha(x)$ goes to zero exponentially, the closed-loop dynamics from $v$ to $y$ is fully linearized despite the fact that $l_y>l_u$, which, in general, is not true in the application of input-output linearizing control to tall systems.  

\begin{figure}[!ht]
    \centering
    \includegraphics[width = 0.6\columnwidth]{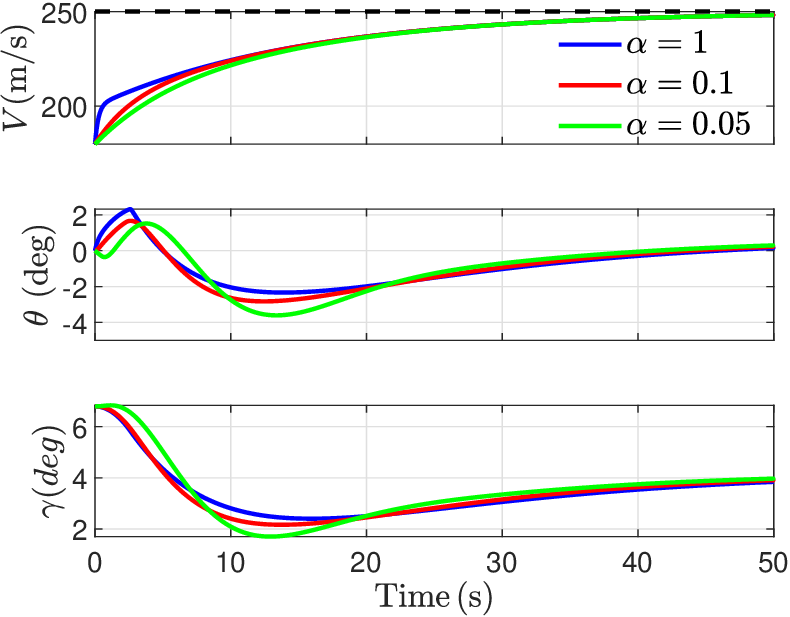}
    \caption{
    Effect of controller gains on the closed-loop response.
    The figure shows the velocity, flight-path angle, and pitch angle response in the case where nominal controller gains are multiplied by a scalar  $\alpha.$ 
    %
    }
    \label{fig:different_gains}
\end{figure}

Next, we scale all controller gains $K_i$ by a scalar factor $\alpha.$
Figure \ref{fig:different_gains} shows the closed-loop response of the aircraft as the IOL controller gains scaled with three different values of $\alpha$.
This example shows that  since all outputs are linearized and decoupled, arbitrary dynamics can be imposed on each output. 

\begin{figure}[!ht]
    \centering
    \includegraphics[width = 0.6\columnwidth]{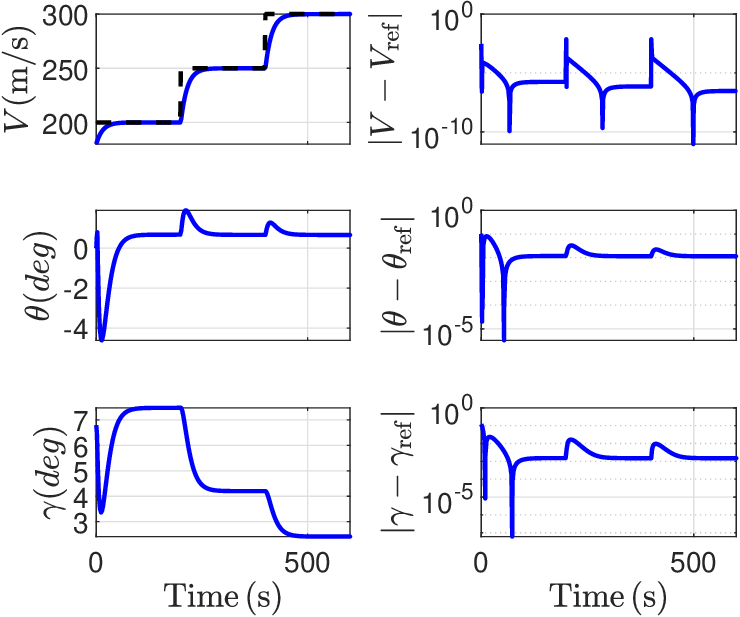}
    \caption{Closed-loop response with biased pitch reference. }
    \label{Disturbed}
\end{figure}

In practice, pitch reference can not be determined exactly. 
We apply the IOL controller in the case where the pitch reference has an unknown bias to model a realistic scenario where the pitch reference is not exactly known.
Such a biased pitch reference may be computed using nominal dynamics. 
Figure \ref{Disturbed} shows the closed-loop response of the aircraft in this case. 
As expected, the output errors do not converge to zero, instead, they converge to nonzero values, which may be sufficient in many practical applications.   

\clearpage
\section{Conclusions and Future Work}
\label{sec:conclusions}
This paper presented an extension of the MIMO input-output linearization method to circumvent the problem of unstable zero dynamics. 
By considering additional outputs, zero dynamics may disappear, but input-output linearization may become impossible if the number of outputs is larger than the number of inputs.
We showed that input-output linearization is still applicable in such tall systems if a geometric condition is satisfied. 
Finally, we showed that in the problem of regulating longitudinal flight dynamics, the additional output removes the zero dynamics, the geometric condition is  satisfied, and thus the longitudinal flight dynamics can be completely linearized.

A current shortcoming of this approach is the requirement of knowledge of various coefficients in the dynamics as well as the exact knowledge of trim conditions. 
In our future work, we will extend the method by including adaptive parameter estimation to estimate the coefficients online.

%

%


\printbibliography

@inproceedings{lambregts1983vertical,
  title={Vertical flight path and speed control autopilot design using total energy principles},
  author={Lambregts, A},
  booktitle={Guidance and Control Conference},
  pages={2239},
  year={1983}
}

@inproceedings{lambregts1983integrated,
  title={Integrated system design for flight and propulsion control using total energy principles},
  author={Lambregts, A},
  booktitle={Aircraft design, systems and technology meeting},
  pages={2561},
  year={1983}
}

@inproceedings{anandakumar2022adaptive,
  title={Adaptive Energy Control of Longitudinal Aircraft Dynamics},
  author={Anandakumar, Ashwin and Bernstein, Dennis and Goel, Ankit},
  booktitle={AIAA SCITECH 2022 Forum},
  pages={0965},
  year={2022}
}

@book{etkin1996dynamics,
  title={Dynamics of flight. Stability and Control},
  author={Etkin, Bernard and Reid, Lloyd D},
  edition = {Third},
  year={1996},
  publisher={John Wiley}
}

@book{stevens2003aircraft,
  title={Aircraft control and simulation: dynamics, controls design, and autonomous systems},
  author={Stevens, Brian L and Lewis, Frank L and Johnson, Eric N},
  edition = {Third},
  year={2003},
  publisher={John Wiley}
}

@article{gavilan2011control,
  title={Control of the longitudinal flight dynamics of an {UAV} using adaptive backstepping},
  author={Gavilan, Francisco and Acosta, JA and Vazquez, Rafael},
  journal={IFAC Proceedings Volumes},
  volume={44},
  number={1},
  pages={1892--1897},
  year={2011},
  publisher={Elsevier}
}

@book{alharbi2019backstepping,
  title={Backstepping Control and Transformation of Multi-Input Multi-Output Affine Nonlinear Systems into a Strict Feedback Form},
  author={Alharbi, Khalid Salim D},
  year={2019},
  publisher={University of Arkansas}
}

@article{kolavennu2001nonlinear,
  title={Nonlinear control of nonsquare multivariable systems},
  author={Kolavennu, S. and Palanki, S. and Cockburn, J. C.},
  journal={Chemical Engineering Science},
  volume={56},
  number={6},
  pages={2103--2110},
  year={2001},
  publisher={Elsevier}
}

@book{isidori1985nonlinear,
  title={Nonlinear control systems: an introduction},
  author={Isidori, Alberto},
  year={1985},
  publisher={Springer}
}

@book{Khalil:1173048,
      author        = "Khalil, Hassan K",
      title         = "{Nonlinear systems; 3rd ed.}",
      publisher     = "Prentice-Hall",
      address       = "Upper Saddle River, NJ",
      year          = "2002",
}

@inproceedings{abdulhamitbilal2007matlab,
  title={Matlab-simulink nonlinear modeling and simulation of aircraft longitudinal dynamics},
  author={Abdulhamitbilal, Erkan and Jafarov, Elbrous M and Kavsao{\u{g}}lu, M {\c{S}}erif},
  booktitle={Proc. of. Eurosim},
  year={2007}
}

\end{document}